\documentclass[a4paper,12pt]{article}
\usepackage{amsmath,amsthm,amsfonts,amssymb}
\usepackage[matrix,arrow]{xy}
\makeatletter
\renewcommand\section{\@startsection {section}{1}{\z@}%
  {-3.5ex \@plus -1ex \@minus -.2ex}%
  {2.3ex \@plus.2ex}%
  {\normalfont\normalsize\scshape}}
\makeatother
\theoremstyle{plain}
\newtheorem{thm}{Theorem}[section]
\newtheorem{prop}[thm]{Proposition}
\newtheorem{lem}[thm]{Lemma}
\newtheorem{cor}[thm]{Corollary}
\theoremstyle{definition}
\newtheorem{dfn}[thm]{Definition}

\newtheorem{ex}[thm]{Example}

\newtheorem{con}[thm]{Conjecture}

\theoremstyle{remark}
\newtheorem{rem}[thm]{Remark}
\renewcommand{\tilde}{\widetilde}
\newcommand{\A}{\mathbb{A}}

\newcommand{\F}{\mathbb{F}}

\renewcommand{\P}{\mathbb{P}}

\newcommand{\R}{\mathbb{R}}
\renewcommand{\S}{\mathbb{S}}
\newcommand{\C}{\mathbb{C}}

\newcommand{\SO}{\mathcal{O}}

\newcommand{\lra}{\longrightarrow}

\renewcommand{\phi}{\varphi}

\newcommand{\PGL}{\mathrm{PGL}}

\newcommand{\Homeo}{\mathrm{Homeo}}
\newcommand{\Diff}{\mathrm{Diff}}

\DeclareMathSymbol{\sdp}{\mathbin}{AMSb}{"6F}

\newcommand{\alg}{\mathrm{alg}}

\title{Rational real algebraic models of\\topological surfaces}
\author{Indranil Biswas and Johannes Huisman}
\begin{document}
\maketitle
\begin{abstract}
  Comessatti proved that the set of all real points of a rational real
  algebraic surface is either a nonorientable surface, or
  diffeomorphic to the sphere or the torus.  Conversely, it is well
  known that each of these surfaces admits at least one rational real
  algebraic model. We prove that they admit exactly one rational real
  algebraic model.  This was known earlier only for the sphere, the
  torus, the real projective plane and the Klein bottle.
\end{abstract}
\begin{quote}\small
  \textit{MSC 2000:} 14P25, 14E07 \hfill\break \textit{Keywords:} Real
  algebraic surface, topological surface, rational surface, rational
  model, birational map, algebraic diffeomorphism, transitivity,
  geometrically rational surface, geometrically rational model
\end{quote}

\section{Introduction}

Let~$X$ be a rational nonsingular projective real algebraic surface.
Then the set~$X(\R)$ of real points of~$X$ is a compact connected
topological surface.  Comessatti showed that $X(\R)$ cannot be an
orientable surface of genus bigger than~$1$. To put it otherwise,
$X(\R)$ is either nonorientable, or it is orientable and diffeomorphic
to the sphere~$S^2$ or the torus~$S^1\times
S^1$~\cite[p.~257]{Comessatti}.

Conversely, each of these topological surfaces admits a \emph{rational
  real algebraic model}, or \emph{rational model} for short.  In other
words, if $S$ is a compact connected topological surface which is
either nonorientable, or orientable and diffeomorphic to the sphere or
the torus, then there is a nonsingular rational projective real
algebraic surface~$X$ such that~$X(\R)$ is diffeomorphic to~$S$.
Indeed, this is clear for the sphere, the torus and the real
projective plane: the real projective surface defined by the affine
equation~$x^2+y^2+z^2=1$ is a rational model of the sphere~$S^2$, the
real algebraic surface~$\P^1\times\P^1$ is a rational model of the
torus~$S^1\times S^1$, and the real projective plane~$\P^2$ is a
rational model of the topological real projective plane~$\P^2(\R)$. If
$S$ is any of the remaining topological surfaces, then $S$ is
diffeomorphic to the $n$-fold connected sum of the real projective
plane, where $n\geq 2$.
A rational model of such a topological surface is the real
surface obtained from~$\P^2$ by blowing up $n-1$ real points.
Therefore, any compact connected topological surface which is either
nonorientable, or orientable and diffeomorphic to the sphere or the
torus, admits at least one rational model.

Now, if $S$ is a compact connected topological surface admitting a
rational model~$X$, then one can construct many other rational models
of~$S$. To see this,
let $P$ and $\overline{P}$ be a pair of complex conjugate
complex points on~$X$. The blow-up~$\tilde{X}$ of~$X$ at $P$~and
$\overline{P}$ is again a rational model of~$S$. Indeed, since~$P$ and
$\overline{P}$ are nonreal points of~$X$, there are open subsets $U$
of $X$ and $V$ of~$\tilde X$ such that
\begin{itemize}
\item $X(\R)\subseteq U(\R)$, $\tilde{X}(\R)\subseteq V(\R)$, and

\item $U$ and $V$ are isomorphic.
\end{itemize}
In particular,
$X(\R)$ and $\tilde X(\R)$ are diffeomorphic. This means that~$\tilde
X$ is a rational model of~$S$
if $X$ is so. Iterating the process, one can
construct many nonisomorphic rational models of~$S$. We
would like to
consider all such models of~$S$ to be equivalent. Therefore, we
introduce the following equivalence relation on the
collection of all rational
models of a topological surface~$S$.

\begin{dfn}\label{def.iso.}
{\rm Let~$X$ and $Y$ be two rational models of a topological 
surface~$S$.
We say that $X$ and~$Y$ are} isomorphic {\rm as rational models of~$S$
if there is a sequence}
$$
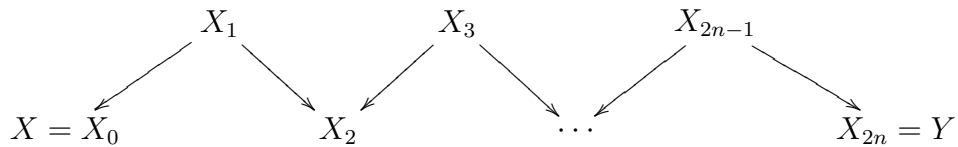
\xymatrix{
&X_1\ar[dl]\ar[dr]&&X_3\ar[dl]\ar[dr]&&X_{2n-1}\ar[dl]\ar[dr]&\\
X=X_0&&X_2&&\cdots&&X_{2n}= Y
}
$$
{\rm where each morphism is a blowing-up at a pair of nonreal complex
conjugate points.}
\end{dfn}

We note that the equivalence relation, in Definition
\ref{def.iso.}, on the collection
of all rational models of a given surface~$S$ is the smallest
one for which the rational
models $X$~and $\tilde X$ mentioned above are equivalent.

Let $X$~and $Y$ be rational models of a topological surface~$S$.
If~$X$ and $Y$ are isomorphic models of~$S$, then the
above sequence of
blowing-ups defines a rational map~
$$
f\colon X\dasharrow Y
$$
having
the following property. There are open subsets $U$ of $X$ and $V$
of~$Y$ such that
\begin{itemize}
\item the restriction of $f$ to $U$ is an isomorphism of
real algebraic varieties from $U$ onto $V$, and

\item $X(\R)\,\subseteq\, U(\R)$ and $Y(\R)\,\subseteq\, V(\R)$.
\end{itemize}
It follows, in particular, that the
restriction of~$f$ to~$X(\R)$ is an \emph{algebraic diffeomorphism}
from~$X(\R)$ onto $Y(\R)$, or in other
words, it is a \emph{biregular map} from $X(\R)$
onto~$Y(\R)$ in the sense of~\cite{BCR}.

Let us recall the notion of an algebraic diffeomorphism.  Let~$X$ and
$Y$ be smooth projective real algebraic varieties.  Then $X(\R)$ and
$Y(\R)$ are compact manifolds, not necessarily connected or nonempty.
Let
\begin{equation}\label{de.f}
f\,\colon\, X(\R)\,\longrightarrow\, Y(\R)
\end{equation}
be a map. Choose affine open subsets $U$ of~$X$ and $V$ of~$Y$ such
that~$X(\R)\subseteq U(\R)$ and $Y(\R)\subseteq V(\R)$.  Since~$U$ and
$V$ are affine, we may assume that they are closed subvarieties
of~$\A^m$ and $\A^n$, respectively.  Then $X(\R)$ is a closed
submanifold of~$\R^m$, and $Y(\R)$ is a closed submanifold of~$\R^n$.
The map~$f$ in~\eqref{de.f} is \emph{algebraic} or \emph{regular} if
there are real polynomials $p_1,\ldots,p_n,q_1,\ldots,q_n$ in the
variables $x_1,\ldots,x_m$ such that none of the polynomials
$q_1,\ldots,q_n$ vanishes on~$X(\R)$, and
$$
f(x)=\left(\frac{p_1(x)}{q_1(x)},\ldots,\frac{p_n(x)}{q_n(x)}\right)
$$
for all~$x\in X(\R)$.

One can check that the algebraicity of~$f$ depends
neither on the choice of the affine open subsets~$U$ and $V$
nor of the choice of the embeddings of $U$~and $V$ in affine space.
Note that the algebraicity of~$f$
immediately implies that~$f$ is a $C^\infty$-map.

The map $f$ in \eqref{de.f} is an \emph{algebraic
diffeomorphism} if $f$ is algebraic, bijective, and $f^{-1}$ is
algebraic.

Again let~$X$ and $Y$ be rational models of a topological surface~$S$.
As observed above, if~$X$ and $Y$ are isomorphic models of~$S$, then
there is an algebraic diffeomorphism
$$
f\colon X(\R)\longrightarrow Y(\R)\, .
$$
Conversely, if there is an algebraic diffeomorphism $f\colon X(\R)
\longrightarrow
Y(\R)$, then $X$ and $Y$ are isomorphic models of~$S$, as it follows
from the well known Weak Factorization Theorem for birational maps
between real algebraic surfaces(see~\cite[Theorem~III.6.3]{BPV} for
the WFT over~$\C$, from which the WFT over $\R$ follows).

Here we address the following question. Given a compact connected
topological surface~$S$, what is the number of nonisomorphic rational
models of~$S$?  

By Comessatti's Theorem, an orientable surface of genus bigger
than~$1$ does not have any rational model.  It is known that the
topological surfaces $S^2$, $S^1\times S^1$ and $\P^2(\R)$ have
exactly one rational model, up to isomorphism (see also
Remark~\ref{rem}).  Mangolte has shown that the same holds for the
Klein bottle~\cite[Theorem~1.3]{Mangolte} (see again
Remark~\ref{rem}).

Mangolte asked how large $n$ should be so that the $n$-fold
connected sum of the real projective plane admits more than
one rational model, up to isomorphism; see the
comments following Theorem 1.3 in \cite{Mangolte}.
The following theorem shows that there is no such integer $n$.

\begin{thm}\label{thmain}
  Let $S$ be a compact connected real two-manifold.
\begin{enumerate}
\item If~$S$ is orientable of genus greater than~$1$, then $S$ does
  not admit any rational model.
\item If~$S$ is either nonorientable, or it is diffeomorphic to one of
  $S^2$ and $S^1\times S^1$, then there is exactly one rational model
  of $S$, up to isomorphism. In other words, any two rational models
  of $S$ are isomorphic.
\end{enumerate}
\end{thm}

Of course, statement~1 is nothing but Comessatti's Theorem referred to
above. Our proof of statement~2 is based on the Minimal Model Program
for real algebraic surfaces developed by J\'anos Koll\'ar in
\cite{Kollar}.  Using this Program, we show that a rational model~$X$
of a nonorientable topological surface~$S$ is obtained from~$\P^2$ by
blowing it up successively in a finite number of real points
(Theorem~\ref{thmmp}).  The next step of the proof of
Theorem~\ref{thmain} involves showing that the model $X$ is isomorphic
to a model~$X'$ obtained from~$\P^2$ by blowing up~$\P^2$ at real
points~$P_1,\ldots,P_n$ of~$\P^2$.  At that point, the proof of
Theorem~\ref{thmain} would have been finished if we were able to prove
that the group~$\Diff_\alg(\P^2(\R))$ of algebraic diffeomorphisms
of~$\P^2(\R)$ acts $n$-transitively on~$\P^2(\R)$.  However, we were
unable to prove such a statement. Nevertheless, a statement we were
able to prove is the following.

\begin{thm}
\label{thntrans}
Let~$n$ be a natural integer.  The group~$\Diff_\alg(S^1\times S^1)$
acts $n$-transitively on~$S^1\times S^1$.
\end{thm}

We conjecture, however, the following.

\begin{con}
  Let~$X$ be a smooth projective rational surface. Let~$n$ be a
  natural integer.  Then the group~$\Diff_\alg(X(\R))$ acts
  $n$-transitively on~$X(\R)$.
\end{con}

The only true evidence we have for the above conjecture is that it holds
for~$X=\P^1\times\P^1$ according to Theorem~\ref{thntrans}.

Now, coming back to the idea of the proof of Theorem~\ref{thmain}, we
know that any rational model of~$S$ is isomorphic to one obtained
from~$\P^2$ by blowing up~$\P^2$ at real points~$P_1,\ldots,P_n$.
Since we have established $n$-transitivity of the group of algebraic
diffeomorphisms of~$S^1\times S^1$, we need to realize~$X'$ as a
blowing-up of~$\P^1\times\P^1$ at a finite number of real points.

Let~$L$ be the real projective line in~$\P^2$ containing $P_1$~and
$P_2$.  Applying a nontrivial algebraic diffeomorphism of~$\P^2$ into
itself, if necessary, we may assume that~$P_i\not\in L$ for~$i\geq3$.
Then we can do the usual transformation of~$\P^2$
into~$\P^1\times\P^1$ by first blowing-up $P_1$~and $P_2$, and then
contracting the strict transform of~$L$.  This realizes~$X'$ as a
surface obtained from~$\P^1\times\P^1$ by blowing-up~$\P^1\times\P^1$
at $n-1$ distinct real points. Theorem~\ref{thmain} then
follows from the $(n-1)$-transitivity of~$\Diff_\alg(S^1\times S^1)$.

We will also address the question of uniqueness of geometrically
rational models of a topological surface. By yet another result of
Comessatti, a geometrically rational real surface~$X$ is rational
if~$X(\R)$ is nonempty and connected. Therefore, Theorem~\ref{thmain}
also holds when one replaces ``rational models'' by ``geometrically
rational models''. Since the set of real points of a geometrically
rational surface is not neccesarily connected, it is natural to study
geometrically rational models of not necessarily connected topological
surfaces. We will show that such a surface has an infinite number of
geometrically rational models, in general.

The paper is organized as follows. In Section~\ref{seHirz} we show
that a real Hirzebruch surface is either isomorphic to the standard
model~$\P^1\times\P^1$ of the real torus~$S^1\times S^1$, or
isomorphic to the standard model of the Klein bottle. The standard
model of the Klein bottle is the real algebraic surface~$B_P(\P^2)$
obtained from the projective plane~$\P^2$ by blowing up one real
point~$P$.  In Section~\ref{semmp}, we use the Minimal Model Program
for real algebraic surfaces in order to prove that any rational model
of any topological surface is obtained by blowing up one of the
following three real algebraic surfaces: $\P^2$, $\S^2$ and
$\P^1\times\P^1$ (Theorem~\ref{thmmp}).  Here~$\S^2$ is the real
algebraic surface defined by the equation~$x^2+y^2+z^2=1$. As a
consequence, we get new proofs of the known facts that the sphere, the
torus, the real projective plane and the Klein bottle admit exactly
one rational model, up to isomorphism of course.  In
Section~\ref{se3fold} we prove a lemma that will have two
applications. Firstly, it allows us to conclude the uniqueness of a
rational model for the ``next'' topological surface, the $3$-fold
connected sum of the real projective plane. Secondly, it also allows
us to conclude that a rational model of a nonorientable topological
surface is isomorphic to a model obtained from $\P^2$ by blowing up a
finite number of distinct real points~$P_1,\ldots,P_n$ of~$\P^2$.  In
Section~\ref{setrans} we prove $n$-transitivity of the group of
algebraic diffeomorphisms of the torus~$S^1\times S^1$.  In
Section~\ref{sealgdiff} we construct a nontrivial algebraic
diffeomorphism~$f$ of~$\P^2(\R)$ such that the real points~$f(P_i)$,
for~$i=3,\ldots,n$, are not on the real projective line
through~$f(P_1)$ and $f(P_2)$. In Section~\ref{semain} we put all the
pieces together and complete the proof of Theorem~\ref{thmain}.  In
Section~\ref{segrm} we show by an example that the uniqueness does not
hold for geometrically rational models of nonconnected topological
surfaces.

\subparagraph{Acknowledgement.} The second author thanks the Tata Institute 
of Fundamental Research for its hospitality.

\section{Real Hirzebruch surfaces}
\label{seHirz}

The set of real points of the rational real algebraic
surface~$\P^1\times\P^1$ is the torus~$S^1\times S^1$. We call this
model the \emph{standard model} of the real torus. Fix a real
point~$O$ of the projective plane~$\P^2$. The rational real algebraic
surface~$B_O(\P^2)$ obtained from~$\P^2$ by blowing up the real
point~$O$ is a model of the Klein bottle~$K$. We call this model the
\emph{standard model} of the Klein bottle.

Let~$d$ be a natural integer.  Let~$\F_d$ be the \emph{real Hirzebruch
  surface} of degree~$d$. Therefore, $\F_d$ is the
compactification~$\P(\SO_{\P^1}(d) \oplus\SO_{\P^1})$ of the line
bundle~$\SO_{\P^1}(d)$ over~$\P^1$. Recall that the real algebraic
surface~$\F_d$ is isomorphic to~$\F_e$ if and only if~$d=e$.  The
restriction of the line bundle~$\SO_{\P^1}(d)$ to the set of real
points~$\P^1(\R)$ of~$\P^1$ is topologically trivial if and only
if~$d$ is even. Consequently, $\F_d$ is a rational model of the
torus~$S^1\times S^1$ if~$d$ is even, and it is a rational model of
the Klein bottle~$K$ if~$d$ is odd (see~\cite[Proposition~VI.1.3]{RAS}
for a different proof).

The following statement is probably well known, and is an easy
consequence of known techniques (compare the proof of Theorem~6.1
in~\cite{Mangolte}).  We have chosen to include the statement and a
proof for two reasons: the statement is used in the proof of
Theorem~\ref{thmmp}, and the idea of the proof turns out also to be
useful in Lemma~\ref{lem}.

\begin{prop}\label{prHirz}
  Let~$d$ be a natural integer.
\begin{enumerate}
\item If~$d$ is even, then $\F_d$ is isomorphic to the standard
  model~$\P^1\times\P^1$ of $S^1\times S^1$.
\item If~$d$ is odd, then $\F_d$ is isomorphic to the standard
  model~$B_O(\P^2)$ of the Klein bottle~$K$.
\end{enumerate}
(All isomorphisms are in the sense of Definition \ref{def.iso.}.)
\end{prop}

\begin{proof}
  Observe that
\begin{itemize}
\item the real algebraic surface~$\P^1\times\P^1$ is
  isomorphic to $\F_0$, and

\item that the real algebraic surface~$B_O(\P^2)$
is isomorphic to~$\F_1$.
\end{itemize}
Therefore, the proposition follows from the following lemma.
\end{proof}

\begin{lem}
  Let $d$ and $e$ be natural integers. Then the
two models $\F_d$ and
  $\F_e$ are isomorphic if and only if~$d\equiv e\pmod2$.
\end{lem}
 
\begin{proof}
  Since the torus is not diffeomorphic to the Klein bottle, the
  rational models $\F_d$ and $\F_e$ are not isomorphic if~$d\not\equiv
  e\pmod2$. Conversely, if~$d\equiv e\pmod2$, then $\F_d$ and $\F_e$
  are isomorphic models, as follows from the following lemma using
  induction.
\end{proof}

\begin{lem} Let~$d$ be a natural integer. The two
rational models $\F_d$ and~$\F_{d+2}$ are isomorphic.
\end{lem}

\begin{proof}
  Let~$E$ be the section at infinity of~$\F_d$.  The self-intersection
  of~$E$ is equal to~$-d$. Choose nonreal complex conjugate points
  $P$~and $\overline{P}$ on~$E$.  Let~$F$~and $\overline{F}$ be the
  fibers of the fibration of~$\F_d$ over~$\P^1$ that contain $P$~and
  $\overline{P}$, respectively. Let~$X$ be the real algebraic surface
  obtained from~$\F_d$ by blowing up $P$~and $\overline{P}$.  Denote
  again by~$E$ the strict transform of~$E$ in~$X$.  The
  self-intersection of~$E$ is equal to~$-d-2$. The strict transforms
  of $F$~and $\overline{F}$, again denoted by $F$~and $\overline{F}$
respectively; they are disjoint smooth rational curves of 
self-intersection~$-1$, and they
  do not intersect~$E$. The real algebraic surface~$Y$ obtained
  from~$X$ by contracting~$F$~and $\overline{F}$ is a
  smooth~$\P^1$-bundle over~$\P^1$.  The image of~$E$ in~$Y$
has self-intersection~$-d-2$.  It follows that~$Y$ is
isomorphic to~$\F_{d+2}$ as a real algebraic surface. Therefore,
we conclude that $\F_d$ and $\F_{d+2}$ are isomorphic models.
\end{proof}

\section{Rational models}
\label{semmp}

Let~$Y$ be a real algebraic surface.  A real algebraic surface~$X$ is
said to be \emph{obtained from~$Y$ by blowing up} if there is a
nonnegative integer~$n$, and a sequence of morphisms
\[
\xymatrix{
X=X_n\ar[r]^{f_n}&X_{n-1}\ar[r]^{f_{n-1}}&\cdots\ar[r]^{f_1}&X_0=Y\, ,
}
\]
such that for each $i=1,\ldots,n$, the morphism $f_i$ is either the
blow up of~$X_{i-1}$ at a real point, or it is the blow up
of~$X_{i-1}$ at a pair of distinct complex conjugate points.

The surface~$X$ is said to be obtained from~$Y$ by blowing up
\emph{at real points only} if for each~$i=1,\ldots,n$,
the morphism~$f_i$ is a blow up of~$X_{i-1}$ at
a real point of~$X_{i-1}$.

One defines,
similarly, the notion of a real algebraic surface obtained from~$Y$ by
blowing up \emph{at nonreal points only}.

The real algebraic surface defined by the affine
equation
$$
x^2+y^2+z^2=1
$$
will be denoted by~$\S^2$. Its set of real
points is the two-sphere~$S^2$. The real Hirzebruch surface~$\F_1$
will be simply denoted by~$\F$.  Its set of real points is the Klein
bottle~$K$.

Thanks to the Minimal Model Program for real algebraic surfaces due to
J\'anos Koll\'ar~\cite[p. 206, Theorem 30]{Kollar}, one has the 
following statement: 

\begin{thm}\label{thmmp}
  Let~$S$ be a compact connected topological surface.  Let~$X$ be a
  rational model of~$S$. 
\begin{enumerate}
\item If $S$ is not orientable then $X$ is isomorphic to a rational
  model of~$S$ obtained from~$\P^2$ by blowing up at real points only.
\item If $S$ is orientable then $X$ is isomorphic to $\S^2$
  or~$\P^1\times\P^1$, as a model.
\end{enumerate}
\end{thm}

\begin{proof}
  Apply the Minimal Model Program to~$X$ in order to obtain a sequence
  of blowing-ups as above, where $Y$ is one of the following:
\begin{enumerate}
\item a minimal surface,

\item a conic bundle
  over a smooth real algebraic curve,

\item a Del Pezzo surface of degree
  $1$ or $2$, and

\item  $\P^2$ or $\S^2$.
\end{enumerate}
(See \cite[p. 206, Theorem 30]{Kollar}.) The surface $X$
being rational,  we
know that $X$ is not bi\-rational to a minimal surface. This rules out
the case of $Y$ being a minimal surface.  Since~$X(\R)$ is connected,
it can be shown that $X$ is not birational to a Del Pezzo surface of
degree $1$~or $2$. Indeed, such Del Pezzo surfaces have disconnected
sets of real points~\cite[p. 207, Theorem 33(D)(c--d)]{Kollar}.
This rules out the
case of~$Y$ being a Del Pezzo surface of degree $1$~or $2$.  It
follows that
\begin{itemize}
\item either $Y$ is a conic bundle, or

\item $Y$ is isomorphic to $\P^2$, or

\item $Y$ is isomorphic to $\S^2$.
\end{itemize}
We will show that the
statement of the theorem holds in all these three cases.

If~$Y$ is isomorphic to~$\P^2$, then $Y(\R)$ is not orientable.  Since
$X$ is obtained from~$Y$ by blowing up,
it follows that $X(\R)$ is not orientable
either. Therefore, the surface
$S$ is not orientable, and also $X$ is isomorphic to a
rational model of~$S$ obtained from~$\P^2$ by blowing up. Moreover, it
is easy to see that $X$ is then isomorphic to a rational model
of~$S$ obtained from~$\P^2$ by blowing up at real points only.  This
settles the case when $Y$ is isomorphic to~$\P^2$.

If~$Y$ is isomorphic to~$\S^2$, then there are two cases to consider:
(1) the case of $S$ being orientable, (2) and the case of $S$ being
nonorientable.  If $S$ is orientable, then $X(\R)$ is orientable too,
and~$X$ is obtained from $Y$ by blowing up at nonreal points only.  It
follows that~$X$ is isomorphic to $\S^2$ as a model.

If $S$ is nonorientable, then $X(\R)$ is nonorientable too, and~$X$ is
obtained from~$\S^2$ by blowing up a nonempty set of real points.
Therefore, the map~$X\longrightarrow Y$ factors through a blow up
$\widetilde{\S}^2$ of $\S^2$ at a real point.  Now, $\widetilde{\S}^2$
contains two smooth disjoint complex conjugated rational curves of
self-intersection~$-1$.  When we contract them, we obtain a real
algebraic surface isomorphic to~$\P^2$.  Therefore, $X$ is obtained
from~$\P^2$ by blowing up.  It follows again that $X$ is isomorphic to
a rational model of~$S$ obtained from~$\P^2$ by blowing up at real
points only.  This settles the case when $Y$ is isomorphic to~$\S^2$.

The final case to consider is the one where $Y$ is a conic bundle over
a smooth real algebraic curve~$B$.  Since $X$ is rational, $B$ is
rational. Moreover, $B$ has real points because $X$ has real points.
Hence, the curve $B$ is isomorphic to $\P^1$.

The singular fibers of the
the conic bundle $Y$ over $B$ are real, and
moreover, the number of singular fibers is even.
Since $X(\R)$ is connected, we conclude that
$Y(\R)$ is connected too. it follows that the conic bundle~$Y$
over~$B$ has either no singular fibers or exactly $2$ singular fibers.
If it has exactly $2$ singular fibers, then~$Y$ is isomorphic
to~$\S^2$ \cite[Lemma~3.2.4]{Kollar2}, a case we have
already dealt with.

Therefore, we may assume that $Y$ is a smooth $\P^1$-bundle
over~$\P^1$. Therefore, $Y$ is a real Hirzebruch surface. By
Proposition~\ref{prHirz}, we may suppose that~$Y=\P^1\times\P^1$, or
that~$Y=\F$. Since~$\F$ is obtained from~$\P^2$ by
blowing up one real point, the case~$Y=\F$ follows from the
case of~$Y=\P^2$ which we have already dealt with above.
 
Therefore, we may assume that~$Y=\P^1\times\P^1$. Again, two cases are
to be considered: (1) the case of~$S$ being orientable, and
(2) the case of  $S$ being nonorientable.
If~$S$ is orientable, $X(\R)$ is orientable, and
$X$ is obtained from~$Y$ by blowing up at non real points only. It
follows that~$X$ is isomorphic as a model to~$\P^1\times\P^1$.  If~$S$
is not orientable, $X$ is obtained from~$Y$ by blowing up, at least,
one real point.  Since~$Y=\P^1\times\P^1$, a blow-up of~$Y$ at one
real point is isomorphic to a blow-up of~$\P^2$ at two real points. We
conclude again by the case of~$Y=\P^2$ dealt with above.
\end{proof}

Note that Theorem~\ref{thmmp} implies Comessatti's Theorem referred to
in the introduction, i.e., the statement to the effect that any
orientable compact connected topological surface of genus greater
than~$1$ does not admit a rational model (Theorem~\ref{thmain}.1).

\begin{rem}\label{rem}
  For sake of completeness let us show how Theorem~\ref{thmmp} implies
  that the surfaces $S^2,S^1\times S^1, \P^2(\R)$ and the Klein bottle
  $K$ admit exactly one rational model. First, this is clear for the
  orientable surfaces~$S^2$ and $S^1\times S^1$.

  Let~$X$ be a rational model of $\P^2(\R)$. From Theorem~\ref{thmmp},
  we know that $X$ is isomorphic to a rational model of~$\P^2(\R)$
  obtained from $\P^2$ by blowing up at real points only.  Therefore,
  we may assume that~$X$ itself is obtained from~$\P^2$ by blowing up
  at real points only.  Since~$X(\R)$ is diffeomorphic to~$\P^2(\R)$,
  it follows that $X$ is isomorphic to~$\P^2$.  Thus any rational
  model of~$\P^2(\R)$ is isomorphic to~$\P^2$ as a model.

  Let~$X$ be a rational model of the Klein bottle~$K$. Using
  Theorem~\ref{thmmp} one may assume that~$X$ is a blowing up of
  $\P^2$ at real points only.  Since~$X(\R)$ is diffeomorphic to the
  $2$-fold connected sum of~$\P^2(\R)$, the surface $X$ is a blowing
  up of $\P^2$ at exactly one real point. It follows that~$X$ is
  isomorphic to~$\F$.  Therefore, any rational model of the Klein
  bottle $K$ is isomorphic to~$\F$, as a model;
  compare with \cite[Theorem~1.3]{Mangolte}.
\end{rem}

One can wonder whether the case where $S$ is a $3$-fold connected sum
of real projective planes can be treated similarly. The first
difficulty is as follows.
It is, a priori, not clear why the following two
rational models of $\#^3\P^2(\R)$ are isomorphic. The first one is
obtained from~$\P^2$ by blowing up two real points of~$\P^2$. The
second one is obtained by a successive blow-up of $\P^2$: first blow
up~$\P^2$ at a real point, and then blow up a real point of the
exceptional divisor. In the next section we prove that these
two models are isomorphic.

\section{The $3$-fold connected sum of the real projective plane}
\label{se3fold}

We start with a lemma.

\begin{lem}\label{lem}
  Let~$P$ be a real point of~$\P^2$, and let~$B_P(\P^2)$ be the
  surface obtained from~$\P^2$ by blowing up~$P$. Let~$E$ be the
  exceptional divisor of~$B_P(\P^2)$ over~$P$.  Let~$L$ be any real
  projective line of~$\P^2$ not containing~$P$. Consider~$L$ as a
  curve in~$B_P(\P^2)$.  Then there is a birational map
  \[
  f\colon B_P(\P^2)\dasharrow B_P(\P^2)
  \]
  whose restriction to the set of real points is an algebraic
  diffeomorphism such that $f(L(\R))=E(\R)$.
\end{lem}

\begin{proof}
  The real algebraic surface~$B_P(\P^2)$ is isomorphic to the real
  Hirzebruch surface~$\F=\F_1$, and any isomorphism between them takes
  the exceptional divisor of~$B_P(\P^2)$ to the section at infinity of
  the conic bundle~$\F/\P^1 = \P(\SO_{\P^1}(1)\oplus\SO_{\P^1})$.
The line~$L$
  in~$B_P(\P^2)$ is given by a unique section of $\SO_{\P^1}(1)$ 
over~$\P^1$; this section of $\SO_{\P^1}(1)$
will also be denoted by~$L$. We denote again by~$E$ the
  section at infinity of~$\F$.

  We have to show that there is a birational self-map~$f$ of~$\F$ such
  that the equality $f(L(\R))=E(\R)$ holds.  Let $R$ be a nonreal
  point of~$L$.  Let~$F$ be the fiber of the conic bundle $\F$ passing
  through~$R$. The blowing-up of~$\F$ at the pair of points $R$~and
  $\overline{R}$ is a real algebraic surface in which we can contract
  the strict transforms of $F$~and $\overline{F}$.  The real algebraic
  surface one obtains after these two contractions is again isomorphic
  to~$\F$.

  Therefore, we have a birational self-map~$f$ of~$\F$ whose
  restriction to the set of real points is an algebraic
  diffeomorphism.  The image, by $f$,
 of the strict transform of~$L$ in~$\F$ has self-intersection
$-1$. Therefore, the image, by $f$, of the strict
  transform of $L$ coincides with $E$. In particular, we have
  $f(L(\R))=E(\R)$.
\end{proof}

\begin{prop}\label{pr3fold}
  Let~$S$ be the $3$-fold connected sum of~$\P^2(\R)$. Then $S$ admits
  exactly $1$ rational model.
\end{prop}

\begin{proof}
  Fix two real points $O_1,O_2$ of~$\P^2$, and let~$B_{O_1,O_2}(\P^2)$
  be the real algebraic surface obtained from~$\P^2$ by blowing up
  $O_1$~and $O_2$. The surface~$B_{O_1,O_2}(\P^2)$ is a rational model
  of the $3$-fold connected sum of~$\P^2(\R)$.
 
  Let~$X$ be a rational model of~$\P^2(\R)$. We prove that~$X$ is
  isomorphic to~$B_{O_1,O_2}(\P^2)$, as a model. By
  Theorem~\ref{thmmp}, we may assume that~$X$ is obtained from~$\P^2$
  by blowing up real points only.  Since~$X(\R)$ is diffeomorphic to a
  $3$-fold connected sum of the real projective plane, the surface $X$
  is obtained from~$\P^2$ by blowing up twice real points.  More
  precisely, there is a real point~$P$ of~$\P^2$ and a real point~$Q$
  of the blow-up~$B_P(\P^2)$ of~$\P^2$ at~$P$, such that~$X$ is
  isomorphic to the blow-up $B_Q(B_P(\P^2))$ of $B_P(\P^2)$ at $Q$.

  Choose any real projective line~$L$ in~$\P^2$ not containing~$P$.
  Then, $L$ is also a real curve in~$B_P(\P^2)$. We may assume
  that~$Q\not\in L$.  By Lemma~\ref{lem}, there is a birational
  map~$f$ from $B_P(\P^2)$ into itself whose restriction to the set of
  real points is an algebraic diffeomorphism, and such
  that
$$
f(L(\R))=E(\R)\, ,
$$
where $E$ is the exceptional divisor
on~$B_P(\P^2)$. Let~$R=f(Q)$. Then~$R\not \in E$, and $f$ induces a
  birational isomorphism
  $$
  \tilde{f}\colon B_Q(B_P(\P^2))\lra B_R(B_P(\P^2))
  $$
  whose restriction to the set of real points is an algebraic
  diffeomorphism. Since~$R\not\in E$,
the point $R$ is a real point of~$\P^2$
  distinct from~$P$, and the blow-up~$B_R(B_P(\P^2))$ is equal to the
  blow up~$B_{P,R}(\P^2)$ of~$\P^2$ at the real points $P,R$
  of~$\P^2$. It is clear that~$B_{P,R}(\P^2)$ is isomorphic
  to~$B_{O_1,O_2}(\P^2)$. It follows that~$X$ is isomorphic
  to~$B_{O_1,O_2}(\P^2)$ as rational models of the $3$-fold connected
  sum of~$\P^2(\R)$.
\end{proof}

\begin{lem}\label{leblow-up}
  Let~$S$ be a nonorientable surface and let~$X$ be a rational model
  of~$S$.  Then, there are distinct real points $P_1,\ldots,P_n$
  of~$\P^2$ such that $X$ is isomorphic to the blowing-up of~$\P^2$ at
  $P_1,\ldots,P_n$, as a model.
\end{lem}

\begin{proof}
  By Theorem~\ref{thmmp}, we may assume that $X$ is obtained
  from~$\P^2$ by blowing up at real points only. Let
\begin{equation}\label{eq1}
\xymatrix{
    X=X_n\ar[r]^{f_n}&X_{n-1}\ar[r]^{f_{n-1}}&\cdots\ar[r]^{f_1}&X_0=\P^2.
  }
\end{equation}
be a sequence of blowing ups, where for each $i=1,\ldots,n$,
the map $f_i$ is a
  blowing up of~$X_{i-1}$ at a real point~$P_i$ of~$X_{i-1}$.

To a sequence of blowing-ups as in \eqref{eq1}
is associated a forest~$F$ of
  trees.  The vertices of~$F$ are the centers~$P_i$ of the blow-ups
  $f_i$.  For~$i>j$, there is an edge between the points $P_i$~and
  $P_j$ in~$F$ if
\begin{itemize}
\item the composition~$f_{j+1}\circ\cdots\circ f_{i-1}$ is an
  isomorphism at a neighborhood of~$P_i$, and

\item maps~$P_i$ to a point
  belonging to the exceptional divisor~$f_j^{-1}(P_j)$ of~$P_j$
  in $X_j$.
\end{itemize}

  Let~$\ell$ be the sum of the lengths of the trees belonging to~$F$.
  We will
show by induction on~$\ell$ that~$X$ is isomorphic, as a model,
  to the blowing-up of~$\P^2$ at a finite number of distinct real
  points of~$\P^2$.

  This is obvious if~$\ell=0$. If~$\ell\neq0$, let~$P_j$ be the root
  of a tree of nonzero length, and let~$P_i$ be the vertex of that
  tree lying immediately above~$P_j$. By changing the order of the
  blowing-ups~$f_i$, we may assume that~$j=1$ and $i=2$.

  Choose a real projective line~$L$ in~$\P^2$ which does not contain
  any of the roots of the trees of~$F$. By Lemma~\ref{lem}, there is a
  birational map~$g_1$ from~$X_1=B_{P_1}(\P^2)$ into itself whose
  restriction to the set of real points is an algebraic diffeomorphism
  and satisfies the condition $g_1(L(\R))=E(\R)$, where $E$ is the
  exceptional divisor of~$X_1$.

  Put~$X_0'=\P^2$, $X_1'=X_1$, and $f_1'=f_1$. We consider~$g_1$ as a
  birational map from~$X_1$ into~$X_1'$.  Put~$P_2'=g_1(P_2)$.
  Let~$X_2'$ be the blowing-up of~$X_1'$ at $P_2'$, and
  let
\[
f_2'\colon X_2'\longrightarrow X_1'
\]
be the blowing-up morphism. Then,
  $g_1$ induces a birational map~$g_2$ from~$X_2$ into $X_2'$ which is
  an algebraic diffeomorphism on the set of real points.

  By iterating this construction, one gets a sequence of blowing ups
\[
f_i'\colon X_i'\longrightarrow  X_{i-1}'\, ,
\]
where $i=1,\ldots,n$, and birational
  morphisms~$g_i$ from $X_i$ into~$X_i'$ whose restrictions to the
  sets of real points are algebraic diffeomorphisms. In particular,
  the rational models $X=X_n$ and $X'=X_n'$ of~$S$ are isomorphic.

  Let~$F'$ be the forest of the trees of centers of~$X'$. Then the sum
  of the lengths~$\ell'$ of the trees of~$F'$ is equal to~$\ell-1$.
  Indeed, one obtains $F'$ from $F$ by replacing the tree~$T$ of~$F$
  rooted at~$P_1$ by the disjoint union of the tree $T\setminus{P_1}$
  and the tree~$\{P_1\}$. This follows from the fact that $P_2'$ does
  not belong to the exceptional divisor of~$f_1'$, and that, no root
  of the other trees of~$F$ belongs to the exceptional divisor
  of~$f_1'$ either.
\end{proof}

As observed in the Introduction, if we are able to prove
the $n$-transitivity of the action of the group~$\Diff_\alg(\P^2(\R))$
on~$\P^2(\R)$, then the statement of Theorem~\ref{thmain} would
follow from Lemma~\ref{leblow-up}.  However, we did not succeed in
proving so. Nevertheless, we will prove the $n$-transitivity
of~$\Diff(S^1\times S^1)$, which is the subject of the next section.

Now that we know that the topological surfaces $S^1,S^1\times S^1$ and
$\#^n\P^2(\R)$, for $n=1,2,3$, admit exactly one rational model, one
may also wonder whether Lemma~\ref{leblow-up} allows
us to tackle the
``next'' surface, which is the $4$-fold connected sum of~$\P^2(\R)$.
We note that
Theorem~\ref{thmain} and Lemma~\ref{leblow-up} imply that a rational
model of such a surface is isomorphic to a surface obtained
from~$\P^2$ by blowing up $3$ distinct real points. However, it it is
not clear why the two surfaces of the following type
are isomorphic as models. Take three
non--collinear real points~$P_1,P_2,P_3$, and three collinear
distinct real points~$Q_1,Q_2,Q_3$ of~$\P^2$. Then the surfaces
$X=B_{P_1,P_2,P_3}(\P^2)$ and $Y=B_{Q_1,Q_2,Q_3}(\P^2)$ are rational
models of~$\#^4\P^2(\R)$ (the $4$-fold connected sum of
$\P^2(\R)$), but it is not clear why they should be
isomorphic. One really seems to need some nontrivial algebraic
diffeomorphism of~$\P^2(\R)$, that maps~$P_i$ to $Q_i$
for~$i=1,2,3$, in order to show that~$X$ and $Y$ are isomorphic
models. We will come back to this in Section~\ref{sealgdiff}
(Lemma~\ref{lealgdiff}).

\section{Algebraic diffeomorphisms of~$S^1\times S^1$ and
  $n$-transitivity}
\label{setrans}

The following statement is a variation on classical polynomial
interpolation.

\begin{lem}\label{leinterpolation}
  Let~$m$ be a positive integer. Let~$x_1,\ldots,x_m$ be distinct real
  numbers, and let~$y_1\ldots,y_m$ be positive real numbers.  Then
  there is a real polynomial $p$ of degree~$2m$ that does not have
  real zeros, and satisfies the condition~$p(x_i)=y_i$ for all~$i$.
\end{lem}

\begin{proof}
Set
\[
p(\zeta):=
\sum_{j=1}^m\prod_{k\neq j}\frac{(\zeta-x_k)^2}{(x_j-x_k)^2}\cdot y_j.
\]
Then~$p$ is of degree~$2m$, and $p$ does not have real zeros.
Furthermore, we have $p(x_i)=y_i$ for all~$i$.
\end{proof}

\begin{cor}\label{cointerpolation}
  Let~$m$ be a positive integer.  Let~$x_1,\ldots,x_m$ be distinct real
  numbers, and let~$y_1\ldots,y_m,z_1,\ldots,z_m$ be positive real
  numbers.  Then there are real polynomials $p$ and $q$ without any
real zeros such that ${\rm degree}(p)= {\rm degree}(q)$, and
\[
\frac{p(x_i)}{q(x_i)}=\frac{y_i}{z_i}
\]
for all $1\leq i\leq m$.\qed
\end{cor}

The interest in the rational functions~$p/q$ of the
above type lies in the following fact.

\begin{lem}\label{leauto}
  Let $p$~and $q$ be two real polynomials of same degree that do not
  have any real zeros. Define the rational
  map~$f\colon\P^1\times\P^1\dasharrow\P^1\times\P^1$ by
  \[
  f(x,y)=\left(x,\frac{p(x)}{q(x)}\cdot y\right).
  \]
  Then~$f$ is a birational map of~$\P^1\times\P^1$ into itself whose
  restriction to the set of real points is an algebraic
  diffeomorphism.\qed
\end{lem}

\begin{thm}\label{thp1p1} 
  Let~$n$ be a natural integer. The group~$\Diff_\alg(\P^1\times\P^1)$
  acts $n$-transitively on~$\P^1(\R)\times\P^1(\R)$.
\end{thm}

\begin{proof}
  Choose $n$ distinct real points~$P_1,\ldots,P_n$ and $n$ distinct
  real points $Q_1,\ldots,Q_n$ of~$\P^1\times\P^1$.  We need to show
  that there is a birational map $f$ from~$\P^1\times\P^1$ into
  itself, whose restriction to~$(\P^1\times\P^1)(\R)$ is an
  algebraic diffeomorphism, such that~$f(P_i)=Q_i$, for $i=1,\ldots,
  n$.

  First of all, we may assume that~$P_1,\ldots,P_n,Q_1,\ldots,Q_n$ are
  contained in the first open quadrant of~$\P^1(\R)\times\P^1(\R)$. In
  other words, the coordinates of $P_i$ and $Q_i$ are strictly
  positive real numbers.  Moreover, it suffices to prove the statement
  for the case where~$Q_i=(i,i)$ for all~$i$.
  
  By the hypothesis above, there are positive real numbers $x_i,y_i$
  such that~$P_i=(x_i,y_i)$ for all~$i$.  By
  Corollary~\ref{cointerpolation}, there are real polynomials~$p$ and
  $q$ without any real zeros such that ${\rm degree}(p)=
{\rm degree}(q)$, and such
  that the real numbers
  \[
  \frac{p(x_i)}{q(x_i)}\cdot y_i
  \]
  are positive and distinct for all~$i$.
  Define~$f\colon\P^1\times\P^1\dasharrow\P^1\times\P^1$ by
  \[
  f(x,y):=\left(x,\frac{p(x)}{q(x)}\cdot y\right).
  \]
  By Lemma~\ref{leauto}, $f$ is birational, and its restriction to
  $(\P^1\times\P^1)(\R)$ is an algebraic diffeomorphism.  By
  construction, the points~$f(P_i)$ have distinct second coordinates.
  Therefore, replacing $P_i$ by~$f(P_i)$ if necessary, we may assume
  that the points~$P_i$ have distinct second coordinates, which
  implies that $y_1,\ldots,y_m$ are distinct positive real numbers.

  By Corollary~\ref{cointerpolation}, there are real polynomials~$p,q$
  without any real zeros such that ${\rm degree}(p)=
{\rm degree}(q)$, and
  \[
  \frac{p(y_i)}{q(y_i)}\cdot x_i=i.
  \]
  Define~$f\colon\P^1\times\P^1\dasharrow\P^1\times\P^1$ by
  \[
  f(x,y)=\left(\frac{p(y)}{q(y)}\cdot x,y\right).
  \]
  By Lemma~\ref{leauto}, $f$ is birational and its restriction to the
  set of real points is an algebraic diffeomorphism. By construction,
  one has~$f(P_i)=i$ for all~$i$.  Therefore, we may assume
  that~$P_i=(i,y_i)$ for all~$i$.

  Now, again by Corollary~\ref{cointerpolation}, there are two real
  polynomials~$p$ and $q$ without any real zeros such that
  $\text{degree}(p)= \text{degree}(q)$, and
  \[
  \frac{p(i)}{q(i)}\cdot y_i=i
  \]
  for all~$i$.  Define~$f\colon\P^1\times\P^1\dasharrow\P^1\times\P^1$
  by
  \[
  f(x,y)=\left(x,\frac{p(x)}{q(x)}\cdot y\right).
  \]
  By Lemma~\ref{leauto}, $f$ is birational, and its restriction
  to the set of real points is an algebraic diffeomorphism. By
  construction~$f(P_i)=Q_i$ for all~$i$.
\end{proof}

\begin{rem}
  One may wonder whether Theorem~\ref{thp1p1} implies that the group
  $\Diff_\alg(\P^2(\R))$ acts $n$-transitively on~$\P^2(\R)$. We will
explain the implication of Theorem~\ref{thp1p1} in that direction.
Let~$P_1,\ldots,P_n$ be
  distinct real points of~$\P^2$, and let~$Q_1,\ldots,Q_n$ be distinct
  real points of~$\P^2$. Choose a real projective line~$L$ in~$\P^2$
  not containing any of the points $P_1,\ldots,P_n,Q_1,\ldots,Q_n$.
  Let~$O_1$ and~$O_2$ be distinct real points of~$L$.
  Identify~$\P^1\times\P^1$ with the surface obtained from~$\P^2$ by,
  first, blowing up~$O_1,O_2$ and, then, contracting the strict
  transform of~$L$. Denote by~$E_1$ and $E_2$ the images of the
  exceptional divisors over~$O_1$~and $O_2$ in~$\P^1\times\P^1$,
  respectively.  We denote again by~$P_1,\ldots,P_n,Q_1,\ldots,Q_n$
  the real points of~$\P^1\times\P^1$ that correspond to the real
  points~$P_1,\ldots,P_n,Q_1,\ldots,Q_n$ of~$\P^2$.

  Now, the construction in the proof of Theorem~\ref{thp1p1} gives
  rise to a birational map~$f$ from~$\P^1\times\P^1$ into itself which
  is an algebraic diffeomorphism on~$(\P^1\times\P^1)(\R)$ and which
  maps~$P_i$ onto~$Q_i$, for $i=1,\ldots,n$. Moreover, if one carries
  out carefully the construction of~$f$, one has
  that~$f(E_1(\R))=E_1(\R)$ and $f(E_2(\R))=E_2(\R)$ and that the real
  intersection point~$O$ of~$E_1$ and $E_2$ in~$\P^1\times\P^1$ is a
  fixed point of~$f$.

  Note that one obtains back $\P^2$ from~$\P^1\times\P^1$ by blowing
  up~$O$ and contracting the strict transforms of $E_1$~and $E_2$.
  Therefore, the birational map~$f$ of~$\P^1\times\P^1$ into itself
  induces a birational map~$g$ of~$\P^2$ into itself. Moreover,
  $g(P_i)=Q_i$. One may think that~$g$ is an algebraic diffeomorphism
  on~$\P^2(\R)$.  However, the restriction of $g$ to the set of real
  points is not necessarily an algebraic diffeomorphism! In fact, $g$
  is an algebraic diffeomorphism on~$\P^2(\R)\setminus\{O_1,O_2\}$.
  The restriction of~$g$ to~$\P^2(\R)\setminus\{O_1,O_2\}$ does admit
  a continuous extension~$\tilde{g}$ to~$\P^2(\R)$, and~$\tilde{g}$ is
  obviously a homeomorphism. One may call $\tilde{g}$ an
  \emph{algebraic homeomorphism}, but $\tilde{g}$ is
\textit{not necessarily
  an algebraic diffeomorphism}. It is not difficult to find explicit
  examples of such algebraic homeomorphisms that are not
  diffeomorphisms.

  That is the reason why we do not claim to have proven
  $n$-transitivity of~$\Diff_\alg(\P^2(\R)$. The only statement
  about~$\P^2(\R)$ the above arguments prove is the $n$-transitivity
  of the group~$\Homeo_\alg(\P^2(\R))$ of algebraic homeomorphisms.
\end{rem}

\section{A nontrivial algebraic diffeomorphism of~$\P^2(\R)$}
\label{sealgdiff}

The nontrivial diffeomorphisms we have in mind are the following.
They have been studied in another recent paper as well~\cite{RV}.

Let~$Q_1,\ldots, Q_6$ be six pairwise distinct complex points
of~$\P^2$ satisfying the following conditions:
\begin{enumerate}
\item the subset $\{Q_1,\cdots ,Q_6\}$ is closed under complex
  conjugation,
\item the subset $\{Q_1,\cdots ,Q_6\}$ does not lie on a complex
  conic,
\item the complex conic passing through any $5$ of these six points is
  nonsingular.
\end{enumerate}
Denote by~$C_1,\ldots,C_6$ the nonsingular complex conics one thus
obtains. These conics are pairwise complex conjugate.  Consider the
real Cremona transformation~$f=f_Q$ of $\P^2$ defined by first
blowing-up~$\P^2$ at $Q_1,\ldots,Q_6$ and then contracting the strict
transforms of~$C_1,\ldots,C_6$. Let~$R_1,\ldots,R_6$ denote the points
of~$\P^2$ that correspond to the contractions of the
conics~$C_1,\ldots,C_6$.

The restriction to~$\P^2(\R)$ of the birational map~$f$ from~$\P^2$
into itself is obviously an algebraic diffeomorphism.

The Cremona transformation~$f$ maps a real projective line, not
containing any of the points~$Q_1,\ldots,Q_6$, to a real rational
quintic curve having $6$ distinct nonreal double points at the
points~$R_1,\ldots,R_6$. Moreover, it maps a real rational quintic
curve in~$\P^2$ having double points at~$Q_1,\ldots,Q_6$ to a real
projective line in~$\P^2$ that does not contain any of the
points~$R_1,\ldots,R_6$.

Observe that the inverse of the Cremona transformation~$f_Q$ is the
Cremona transformation~$f_R$. It follows that~$f=f_Q$ induces a
bijection from the set of real rational quintics having double
points at~$Q_1,\ldots,Q_6$ onto the set of real projective lines
in~$\P^2$ that do not contain any of~$R_1,\ldots,R_6$.

This section is devoted to the proof of following lemma.

\begin{lem}\label{lealgdiff}
  Let~$n$ be a natural integer bigger than~$1$.  Let~$P_1,\ldots,P_n$
  be distinct real points of~$\P^2$. Then there is a birational map
  of~$\P^2$ into itself, whose restriction to the set of real points
  is an algebraic diffeomorphism, such that the image
  points~$f(P_3),\ldots,f(P_n)$ are not contained in the real
  projective line through $f(P_1)$~and $f(P_2)$.
\end{lem}

\begin{proof}
  Choose complex points~$Q_1,\ldots,Q_6$ of~$\P^2$ as above.  As
  observed before, the Cremona transformation~$f=f_Q$ induces a
  bijection from the set of real rational quintic curves having double
  points at~$Q_1,\ldots,Q_6$ onto the set of real projective lines
  of~$\P^2$ not containing any of the above points~$R_1,\ldots,R_n$.
  In particular, there is a real rational quintic curve~$C$ in~$\P^2$
  having $6$ nonreal double points at~$Q_1,\ldots,Q_6$.

  We show that there is a real projectively linear
  transformation~$\alpha$ of~$\P^2$ such that~$\alpha(C)$ contains
  $P_1$~and $P_2$, and does not contain any of the
  points $P_3,\ldots,P_n$.  The Cremona transformation~$f_{\alpha(Q)}$
  will then be a birational map of~$\P^2$ into itself that has the
  required properties.

  First of all, let us prove that there is~$\alpha\in\PGL_3(\R)$ such
  that~$P_1,P_2\in\alpha(C)$. This is easy. Since~$C$ is a quintic
  curve, $C(\R)$ is infinite. In particular, $C$ contains two distinct
  real points. It follows that there is~$\alpha\in\PGL_3(\R)$ such
  that~$P_1,P_2\in\alpha(C)$.  Replacing $C$ by~$\alpha(C)$ if
  necessary, we may suppose that~$P_1,P_2\in C$.

  We need to show that there is~$\alpha\in\PGL_3(\R)$ such that
  $\alpha(P_1)=P_1$, $\alpha(P_2)=P_2$ and~$\alpha(C)$ does not
  contain any of the points~$P_3,\ldots,P_n$.

  To prove the existence of $\alpha$ by contradiction, assume that
  there is no such automorphism of $\P^2$. Therefore, for all
  $\alpha\in\PGL_3(\R)$ having $P_1$ and $P_2$ as fixed points, the
  image $\alpha(C)$ contains at least one of the points
  of~$P_3,\ldots,P_n$.  Let~$G$ be the stabilizer of the
  pair~$(P_1,P_2)$ for the diagonal action of~$\PGL_3$
  on~$\P^2\times\P^2$.  It is easy to see that $G$ is a geometrically
  irreducible real algebraic group.  Let
  \[
  \rho\colon C\times G\lra\P^2
  \]
  be the morphism defined by~$\rho(P,\alpha)=\alpha(P)$.  Let
$$
X_i := \rho^{-1}(P_i)
$$
be the inverse image, where $i=3,\ldots,n$.
  Therefore, $X_i$ is a real algebraic subvariety of~$C\times G$.  By
  hypothesis, for every $\alpha\in G(\R)$, there is an integer~$i$
  such that $\alpha(C)$ contains~$P_i$.  Denoting by $p$ the
  projection on the second factor from~$C\times G$ onto~$G$, this
  means that
  \[
  \bigcup_{i=3}^n p(X_i(\R))=G(\R).
  \]

  Since~$G(\R)$ is irreducible, there is an integer~$i_0\in [3, n]$ such
  that the semi-algebraic subset~$p(X_{i_0}(\R))$ is Zariski dense
  in~$G(\R)$.  Since~$G$ is irreducible and~$p$ is proper, one
  has~$p(X_{i_0})=G$.  In particular, $P_{i_0}\in \alpha(C)$ for
  all~$\alpha\in G(\C)$. To put it otherwise, $\alpha^{-1}(P_{i_0})\in
  C$ for all~$\alpha\in G(\C)$, which means that the orbit
  of~$P_{i_0}$ under the action of~$G$ is contained in~$C$.  In
  particular, the dimension of the orbit of~$P_{i_0}$ is at most one.
  It follows that~$P_1,P_2$ and $P_{i_0}$ are collinear.  Let~$L$ be
  the projective line through~$P_1,P_2$. Then the orbit of~$P_{i_0}$
  coincides with $L\setminus\{P_1,P_2\}$.  It now follows
  that~$L\subseteq C$. This is in contradiction with the fact that $C$
  is irreducible.
\end{proof}

\section{Proof of Theorem \ref{thmain}.2}
\label{semain}

Let~$S$ be a topological surface, either nonorientable or of genus
less than~$2$.  We need to show that any two rational models of~$S$
are isomorphic. By Remark~\ref{rem}, we may assume that~$S$ is the
$n$-fold connected sum of~$\P^2(\R)$, where~$n\geq3$.

Let~$O_1,\ldots,O_{n-2}$ be fixed pairwise distinct real points
of~$\P^1\times\P^1$, and let~$B_{n-2}(\P^1\times\P^1)$ be the surface
obtained from~$\P^1\times\P^1$ by blowing up the
points~$O_1,\ldots,O_{n-2}$.  It is clear
that~$B_{n-2}(\P^1\times\P^1)$ is a rational model of~$S$.

Now, it suffices to show that any rational model of~$S$ is isomorphic
to~$B_{n-2}(\P^1\times\P^1)$, as a model. Let~$X$ be any rational
model of~$S$.  By Lemma~\ref{leblow-up}, we may assume that there are
distinct real points~$P_1,\ldots,P_m$ of~$\P^2$ such that~$X$ is the
surface obtained from~$\P^2$ by blowing up~$P_1,\ldots,P_m$. Since~$X$
is a rational model of an $n$-fold connected sum of~$\P^2(\R)$, one
has~$m=n-1$. In particular, $m\geq2$. By Lemma~\ref{lealgdiff}, we may
assume that the points~$P_3,\ldots,P_m$ are not contained in the real
projective line~$L$ through $P_1$~and~$P_2$.

The blow-up morphism~$X\longrightarrow \P^2$ factors through
the blow up~$\widetilde{\P}^2=B_{P_1,P_2}(\P^2)$.  The strict
transform~$\widetilde{L}$ of~$L$ has self-intersection $-1$
in~$\widetilde{\P}^2$. If we contract~$\widetilde{L}$, then we obtain
a surface isomorphic to~$\P^1\times\P^1$. Therefore, $X$ is isomorphic
to a model obtained from~$\P^1\times\P^1$ by blowing up $m-1=n-2$
distinct real points of~$\P^1\times\P^1$.  It follows from
Theorem~\ref{thp1p1} that~$X$ is isomorphic
to~$B_{n-2}(\P^1\times\P^1)$.\qed

\section{Geometrically rational models}
\label{segrm}

Recall that a nonsingular projective real algebraic surface~$X$ is
\emph{geometrically rational} if the complex surface~$X_\C
= X\times_{\mathbb R}{\mathbb C}$ is
rational.  Comessatti showed that, if~$X$ is a geometrically rational
real algebraic surface with~$X(\R)$ connected, then~$X$ is rational;
see Theorem~IV of~\cite{Comessatti1} and the remarks thereafter
(see also \cite[Corollary~VI.6.5]{RAS}). Therefore, the main result,
namely
Theorem~\ref{thmain}, also applies to geometrically rational models.
More precisely, we have the following consequence.

\begin{cor}
  Let $S$ be a compact connected real two-manifold. 
\begin{enumerate}
\item If~$S$ is orientable and the genus of~$S$ is greater than~$1$,
  then $S$ does not admit a geometrically rational real algebraic model.
\item If~$S$ is either nonorientable, or it is diffeomorphic to one of
  $S^2$ and $S^1\times S^1$, then there is exactly one geometrically
  rational model of $S$, up to isomorphism. In other words, any two
  geometrically rational models of $S$ are isomorphic.\qed
\end{enumerate}
\end{cor}

Now, the interesting aspect about geometrically rational real surfaces
is that their set of real points can have an arbitrary number of
connected components. More precisely, Comessati proved the following
statement~\cite[p.~263 and further]{Comessatti} (see
also~\cite[Proposition~VI.6.1]{RAS}).

\begin{thm}
  Let~$X$ be a geometrically rational real algebraic surface such
  that~$X(\R)$ is not connected. Then each connected component
  of~$X(\R)$ is either nonorientable or diffeomorphic to~$S^2$.
  Conversely, if~$S$ is a nonconnected compact topological surface
  each of whose connected components is either nonorientable or
  diffeomorphic to~$S^2$, then there is a geometrically rational real
  algebraic surface~$X$ such that~$X(\R)$ is diffeomorphic to~$S$.\qed
\end{thm}

Let~$S$ be a nonconnected topological surface. One may wonder whether
the geometrically rational model of~$S$ whose existence is claimed
above, is unique up to isomorphism of models. The answer is negative,
as shown by the following example.

\begin{ex}
  Let~$S$ be the disjoint union of a real projective plane and $4$
  copies of~$S^2$. Then, any minimal real Del Pezzo surface of
  degree~$1$ is a geometrically rational model
  of~$S$~\cite[Theorem~2.2(D)]{Kollar2}. Minimal real Del Pezzo
  surfaces of degree~$1$ are rigid; this means that any birational map
  between two minimal real Del Pezzo surfaces of degree~$1$ is an
  isomorphism of real algebraic
  surfaces~\cite[Theorem~1.6]{Iskovskikh}. Now, the set of isomorphism
  classes of minimal real Del Pezzo surfaces of degree~$1$ is in
  one-to-one correspondence with the quotient
  set~$\P^2(\R)^8/\PGL_3(\R)$ for the diagonal action of the
  group~$\PGL_3(\R)$. It follows that the topological
  surface~$S$ admits a
  $8$-dimensional continuous family of nonisomorphic geometrically
  rational models. In particular, the number of nonisomorphic
  geometrically rational models of~$S$ is infinite.
\end{ex}

\par\bigskip\noindent\footnotesize
\textsc{School of Mathematics, Tata Institute of Fundamental
Research, Homi Bhabha Road, Mumbai 400005, India\hfill\break
Email:}
\texttt{indranil@math.tifr.res.in}
\medskip
\par\noindent\textsc{D\'epartement de Math\'ematiques,
Laboratoire CNRS UMR 6205,
Universit\'e de Bretagne Occidentale,
6 avenue Victor Le Gorgeu,
CS 93837,
29238 Brest cedex 3,
France\hfil\break E-mail:}
\texttt{johannes.huisman@univ-brest.fr}\hfil\break \textsc{Home page:}
\texttt{http://stockage.univ-brest.fr/$\sim$huisman/}

\par\bigskip\noindent\small
Typeset by \AmS-\LaTeX\
and \Xy-pic
\end{document}